# Technical Report

Privacy-Preserving Economic Dispatch in Competitive Electricity Market

**Lei Wu, Clarkson University**
December 7th 2015

This report has been conducted and completed by Dr. Lei Wu on **December 7th 2015** at Clarkson University, Potsdam NY.

This report includes some initial thoughts on how to use **transformation based linear programming** approach to solve **economic dispatch problems** in competitive electricity markets, in order to disguise critical financial data of individual market participants and physical information of the electricity grid from being disclosed publically, especially when system operators are relying more and more on **cloud computing services** to perform critical and time-consuming power system optimization tasks.

A good portion of this report will be disseminated via a conference paper, and further work on using the transformation based linear programming approach to solve other critical power system operation problems has been presented in another manuscript that is under review by IEEE Transactions on Smart Grid.

# Privacy-Preserving Economic Dispatch in Competitive Electricity Market

Lei Wu, *Senior Member, IEEE*

*Abstract*— With the emerging of smart grid techniques, cyber attackers may be able to gain access to critical energy infrastructure data and strategic market participants may be able to identify their rival producers' offer prices. This paper discusses a privacy-preserving economic dispatch (ED) approach in competitive electricity market, in which individual generation companies (GENCOs) and load serving entities (LSEs) can mask their actual bidding information and physical data by multiplying with random numbers before submitting to Independent System Operators (ISOs)/Regional Transmission Owners (RTOs). This would avoid potential information leakage of critical energy infrastructure and financial data of market participants. The optimal solution to the original ED problem, including optimal dispatches of generators/loads and locational marginal prices (LMPs), can be retrieved from the optimal solution of the proposed privacy-preserving ED approach. Numerical case studies show the effectiveness of the proposed approach for protecting private information of individual market participants while guaranteeing the same optimal ED solution. Computation and communication costs of the proposed privacy-preserving ED approach and the original ED are also compared in case studies.

## NOMENCLATURE

**Variables:**

| | |
|---|---|
| $a, b$ | Indices of buses |
| $D_{vkt}$ | Dispatch of load $v$ at segment $k$ of hour $t$ |
| $i, j$ | Index of GENCOs/LSEs |
| $k, l, t$ | Index of segments/lines/time periods |
| $P_{ukt}$ | Dispatch of unit $u$ at segment $k$ of hour $t$ |
| $u, v$ | Index of generators/loads |
| $\theta_{at}, \theta_{bt}$ | Bus angles of buses $a$ and $b$ at hour $t$ |

**Constants:**

| | |
|---|---|
| $B, L, T$ | Number of buses/lines/time periods |
| $c_{uk}$ | Bidding price of generator $u$ at segment $k$ |
| $d_{vk}$ | Bidding price of load $v$ at segment $k$ |
| $\underline{D}_{vk}, \overline{D}_{vk}$ | Min/max bidding quantity of load $v$ at segment $k$ |
| $G, D$ | Number of GENCOs/LSEs |
| $m_{i(j)}/n_{i(j)}$ | Numbers of constraints/variables related to GENCO $i$ (LSE $j$) |
| $\underline{P}_{uk}, \overline{P}_{uk}$ | Min/max bidding quantity of unit $u$ at segment $k$ |
| $P_u^{dn}, P_u^{up}$ | Ramping down/up limits of unit $u$ |
| $\overline{PL}_l$ | Capacity of line $l$ |
| $x_{ab}$ | Reactance of a line that connects buses $a$ and $b$ |

**Vectors and Matrices:**

| | |
|---|---|
| $\mathbf{B}_b$ | Row of nodal admittance matrix related to bus $b$ |
| $\mathbf{B}, \boldsymbol{\theta}$ | Nodal admittance matrix and bus angle vector |
| $\mathbf{c}_i, \mathbf{d}_j$ | Vector of bidding prices for GENCO $i$/LSE $j$ |
| $\mathbf{D}_j$ | Dispatch vector of loads in LSE $j$ |
| $\mathbf{E}_i, \mathbf{M}_i$ | Constraint coefficients of GENCO $i$ |
| $\mathbf{F}_j, \mathbf{N}_j$ | Constraint coefficients of LSE $j$ |
| $\mathbf{G}$ | Branch susceptance matrix |
| $\mathbf{KL}$ | Branch-bus incidence matrix |
| $\mathbf{KP}_i$ | Bus-generation incidence matrix of GENCO $i$ |
| $\mathbf{KD}_j$ | Bus-load incidence matrix of LSE $j$ |
| $\mathbf{P}_i$ | Dispatch vector of generators in GENCO $i$ |
| $\overline{\mathbf{PL}}$ | Vector of line capacity limits |
| $\mathbf{R}_{l1}, \mathbf{R}_{l2}$ | $(T \cdot L) \times (T \cdot L)$ random positive diagonal matrices generated by ISO |
| $\mathbf{RD}_j, \mathbf{RG}_i$ | An $m_j \times m_j / m_i \times m_i$ random positive diagonal matrix generated by LSE $j$/GENCO $i$ |
| $\mathbf{X}_{l1}, \mathbf{X}_{l2}$ | $(T \cdot L) \times (T \cdot L)$ random positive matrices generated by ISO |
| $\mathbf{X}_b$ | An $(T \cdot B) \times (T \cdot B)$ random positive matrix generated by ISO |
| $\mathbf{XD}_j, \mathbf{XG}_i$ | An $m_j \times m_j / m_i \times m_i$ random matrix generated by LSE $j$/ GENCO $i$ |
| $\mathbf{YD}_j, \mathbf{YG}_i$ | An $n_j \times n_j$ / $n_i \times n_i$ random positive matrix generated by LSE $j$/ GENCO $i$ |
| $\mathbf{Y}_\theta$ | An $n_{ISO} \times n_{ISO}$ random positive matrix generated by ISO |
| $\boldsymbol{\theta}_t$ | Bus angle vector at hour $t$ |
| $\mathbf{0}$ | A zero vector or matrix |
| ' | Transpose of a vector or matrix |
| $\boldsymbol{\alpha}, \boldsymbol{\beta}, \boldsymbol{\mu}, \boldsymbol{\lambda}$ | Vectors of dual variables |
| $\widetilde{sG}, \widetilde{sD}, \widetilde{sL}$ | Vectors of slack variables |

## I. INTRODUCTION

IN restructured power systems, ISOs/RTOs coordinate competitive market participants for ensuring secure system operation and economic market operation. Such operation decisions are usually made via optimization approaches, such as unit commitment (UC) and ED problems. This paper focuses on the ED problem, which determines the least-cost operation of power systems by dispatching generation resources to supply system loads, while satisfying prevailing system-level and unit constraints.

Emerging smart grid techniques have been promoting the increasing participation of new entities into power markets, including smalls-scale renewable/gas-fired generation owners, prosumers equipped with self-owned generation resources, and active customers with demand response capabilities. For instance, NYISO allows various entities to participate in its energy, ancillary service, and demand response markets with the minimum asset size of 1MW [1]. The participation of small-scale entities would significantly increase the volume of data and the scale of optimization models, which boosts the request on high-performance efficient computing architectures. In turn, ISOs/RTOs have been looking into the solution of cloud computing technology for solving their ever increasing operation problems [2]. However, such optimization models always contain critical energy infrastructure information and private financial information of individual market participants,

L. Wu is with Electrical and Computer Engineering Department, Clarkson University, Potsdam, NY 13699 USA. (E-mail: lwu@clarkson.edu).

and in turn security concerns area raised when adopting could computing. Different security schemes, including identity and access management, security group control, and secure data transmission, are being explored regarding the appropriate information disclosure under cloud computing.

Under the restructured power market environment, in order to solve the ED problem, detailed information of individual market participants needs to be collected by the ISO/RTO for the centralized optimization. These market participants are distributed in the system and tend to autonomously maximize their own profits. They share system-level global constraints as well as local constraints that refer to its private limitations or capacities, and serious privacy problems may arise if such information is revealed. Although ISOs/RTOs are neutral and have strict policies regarding the release of such information, there are at least two potential ways that may lead to information leakage: cyberattack and other strategic ways to derive such information. Smart grid technology makes the system even more vulnerable, and attackers could gain access to critical data streams via either hacking into the communication infrastructure, ISOs/RTOs' database, or the cloud. Another way is that a strategic market participant may be able to identifying its rival producers' offer prices with limited available information. This in turn would result in great loss of possible confidential information leakage. For instance, [3] adopted the inverse optimization to reveal price offers of rival producers, with the assumption that accepted generation/ demand blocks for producers/consumers and LMPs are known.

Over the past few years, various schemes have been explored for solving ED problems in a distributed fashion, which does not require a subsystem to disclose its confidential financial information to other subsystems or ISOs/RTOs. Lagrangian Relaxation (LR) was used in [4]-[6] to relax coupling constraints among different subsystems and allow subproblems to be solved separately. Auxiliary problem principle (APP) was applied on the distributed optimal power flow problem [7]-[8], which solves a sequence of auxiliary problems involving the augmented LR. Alternating direction method of multipliers (ADMM) was also utilized for solving distributed convex power system optimization problems [9]. Another important method applied to distributed optimization is the consensus algorithm. Distributed ED approach based on the incremental cost (IC) consensus was discussed in [10]-[11].

Inspired by [12]-[13], different from distributed optimization, this paper proposes a privacy-preserving ED approach which is still a centralized problem solved by ISOs/RTOs. Individual market participants can mask their actual bidding information and physical data by multiplying with random numbers before submitting to ISOs/RTOs. It also allows the ISO to encrypt the physical transmission network information, which makes it possible to rely on cloud computing service providers for performing critical and time-consuming power system optimization tasks while not worry about critical information leakage. In addition, the proposed privacy-preserving ED approach can retrieve the optimal solution to the original ED problem, including optimal dispatches of generators/loads and LMPs. This allows ISOs to continue performing market operation functionalities. Numerical case studies validate the effectiveness of the proposed approach for protecting private information while guaranteeing the same optimal ED solution. Computation and communication costs of the proposed privacy-preserving ED approach and the original ED are also compared in case studies to demonstrate the scalability of the proposed approach. Although this paper focuses on the ED problem for a single ISO/RTO, the proposed approach can also be applied for solving multi-area coordinated ED problems for multiple interconnected ISOs/RTOs [14].

The rest of this paper is organized as follows. Section II presents the proposed privacy preserving ED problem approach. Numerical case studies are presented in Section III, and Section IV summarizes conclusions.

## II. PRIVACY-PRESERVING ED APPROACH

### A. Original ED Problem

The original ED problem is formulated as (1)-(6). The objective (1) is to maximize the social welfare. (2)-(3) represents capacity limits and ramp up/down constraints of individual generators. (4) defines limits for individual loads. (5) enforces transmission line capacity limits. (6) represents the nodal power balance for the DC power flow calculation. Locational marginal prices (LMPs) are determined via dual variables of constraint (6).

$$max \left( -\sum_t \sum_u \sum_k c_{uk} \cdot P_{ukt} + \sum_t \sum_v \sum_k d_{vk} \cdot D_{vkt} \right) \quad (1)$$

$$\underline{P}_{uk} \leq P_{ukt} \leq \bar{P}_k \quad (2)$$

$$P_u^{dn} \leq \sum_k P_{ukt} - \sum_k P_{uk(t-1)} \leq P_u^{up} \quad (3)$$

$$\underline{D}_{vk} \leq D_{vkt} \leq \bar{D}_{vk} \quad (4)$$

$$-\overline{PL}_l \leq (\theta_{at} - \theta_{bt})/x_{ab} \leq \overline{PL}_l \quad (5)$$

$$\sum_{u \in b} \sum_k P_{ukt} - \sum_{v \in b} \sum_k D_{vkt} = \mathbf{B}_b \cdot \boldsymbol{\theta}_t \quad (6)$$

In order to solve the ED problem (1)-(6) in a centralized way, all information from GENCOs and LSEs are to be sent to the ISO/RTO. Such information includes bidding price $c_{uk}$ and physical parameters $\underline{P}_{uk}$, $\bar{P}_{uk}$, $P_u^{dn}$, and $P_u^{up}$ of individual generators, as well as bidding price $d_{vk}$ and physical parameters $\underline{D}_{vk}$ and $\bar{D}_{vk}$ of individual loads. Transmitting such information to ISOs/RTOs and storing them at the ISO/RTO side may be risky to information leakage. Moreover, if the ISO adopts could computing for the ED problem, additional security concerns regarding the critical energy infrastructure information leakage, including $x_{ab}$, $\overline{PL}$, and $\mathbf{B}_b$, may raise.

### B. Privacy-Preserving Linear Programming (LP) Problem

Two privacy-preserving transformation approaches for the LP problem have been investigated in literature [15]-[16], in which certain information owned by individual entities is encrypted via the random matrix transformation before being disclosed. After solving the transformed LP problem, each entity can decode the corresponding solution components and derive the exact solution of the original LP problem. Without loss of generality, the LP problem (7)-(8) is used for the ease of discussion in this subsection. Assuming that boundaries of variables are included in the constraint set (8).

$$min \; [c_{11} \; c_{12} \; c_{21} \; c_{22} \; c_{31} \; c_{32}] \cdot [x_{11} \; x_{12} \; x_{21} \; x_{22} \; x_{31} \; x_{32}]' \quad (7)$$

$$\begin{bmatrix} a_{11} & a_{12} & 0 & 0 & 0 & 0 \\ a_{13} & a_{14} & 0 & 0 & 0 & 0 \\ 0 & 0 & a_{21} & a_{22} & 0 & 0 \\ 0 & 0 & a_{23} & a_{24} & 0 & 0 \\ a_{15} & a_{16} & a_{25} & a_{26} & a_{31} & a_{32} \\ a_{17} & a_{18} & a_{27} & a_{28} & a_{33} & a_{34} \end{bmatrix} \cdot \begin{bmatrix} x_{11} \\ x_{12} \\ x_{21} \\ x_{22} \\ x_{31} \\ x_{32} \end{bmatrix} \begin{matrix} = b_{11} \\ \le b_{12} \\ = b_{21} \\ \le b_{22} \\ = b_{31} \\ \le b_{32} \end{matrix} \quad (8)$$

### B.1. Vertically Partitioned LP Problem

For the LP problem (7)-(8), assuming $c_{11}$, $c_{12}$, and $a_{11}$-$a_{18}$ are owned by Entity I, $c_{21}$, $c_{22}$, and $a_{21}$-$a_{28}$ are owned by Entity II, and $c_{31}$, $c_{32}$, and $a_{31}$-$a_{34}$ are owned by Entity III. Each entity is unwilling to make the corresponding data public to avoid potential private information leakage. This is called a vertically partitioned LP problem, as parameters in each column of the objective and constraints may be owned by different entities.

Entities I-III can generate *2x2* random vectors $\mathbf{Y}_1$, $\mathbf{Y}_2$, and $\mathbf{Y}_3$, respectively. Here, "2" represents the number of variables owned by each entity. That is, $\mathbf{Y}_i = \begin{bmatrix} y_{i1} & y_{i3} \\ y_{i2} & y_{i4} \end{bmatrix}$ for $i = 1,2,3$. Thus, instead of disclosing $c_{11}$, $c_{12}$, and $a_{11}$-$a_{18}$, Entity I provides $[c_{11}\ c_{12}] \cdot \mathbf{Y}_1$ and $\begin{bmatrix} a_{11} & a_{13} & 0 & 0 & a_{15} & a_{17} \\ a_{12} & a_{14} & 0 & 0 & a_{16} & a_{18} \end{bmatrix}' \cdot \mathbf{Y}_1$. Similarly, Entity II provides $[c_{21}\ c_{22}] \cdot \mathbf{Y}_2$ and $\begin{bmatrix} 0 & 0 & a_{21} & a_{23} & a_{25} & a_{27} \\ 0 & 0 & a_{22} & a_{24} & a_{26} & a_{28} \end{bmatrix}' \cdot \mathbf{Y}_2$ instead of $c_{21}$, $c_{22}$, and $a_{21}$-$a_{28}$, and Entity III provides $[c_{31}\ c_{32}] \cdot \mathbf{Y}_3$ and $\begin{bmatrix} 0 & 0 & 0 & 0 & a_{31} & a_{33} \\ 0 & 0 & 0 & 0 & a_{32} & a_{34} \end{bmatrix}' \cdot \mathbf{Y}_3$ instead of $c_{31}$, $c_{32}$, and $a_{31}$-$a_{34}$.

With the encrypted information, the original LP problem (7)-(8) is formulated as a transformed LP problem (9)-(10).

$$min\ [c_{11}\ c_{12}\ c_{21}\ c_{22}\ c_{31}\ c_{32}] \cdot \begin{bmatrix} \mathbf{Y}_1 & 0 & 0 \\ 0 & \mathbf{Y}_2 & 0 \\ 0 & 0 & \mathbf{Y}_3 \end{bmatrix} \cdot$$
$$[\tilde{x}_{11}\ \tilde{x}_{12}\ \tilde{x}_{21}\ \tilde{x}_{22}\ \tilde{x}_{31}\ \tilde{x}_{32}]' \quad (9)$$

$$\begin{bmatrix} a_{11} & a_{12} & 0 & 0 & 0 & 0 \\ a_{13} & a_{14} & 0 & 0 & 0 & 0 \\ 0 & 0 & a_{21} & a_{22} & 0 & 0 \\ 0 & 0 & a_{23} & a_{24} & 0 & 0 \\ a_{15} & a_{16} & a_{25} & a_{26} & a_{31} & a_{32} \\ a_{17} & a_{18} & a_{27} & a_{28} & a_{33} & a_{34} \end{bmatrix} \cdot \begin{bmatrix} \mathbf{Y}_1 & 0 & 0 \\ 0 & \mathbf{Y}_2 & 0 \\ 0 & 0 & \mathbf{Y}_3 \end{bmatrix} \cdot \begin{bmatrix} \tilde{x}_{11} \\ \tilde{x}_{12} \\ \tilde{x}_{21} \\ \tilde{x}_{22} \\ \tilde{x}_{31} \\ \tilde{x}_{32} \end{bmatrix} \begin{matrix} = b_{11} \\ \le b_{12} \\ = b_{21} \\ \le b_{22} \\ = b_{31} \\ \le b_{32} \end{matrix} \quad (10)$$

Comparing (7)-(8) with (9)-(10), it is clear that the extrema are equal. In addition, with the optimal solution $[\tilde{x}_{11}\ \tilde{x}_{12}\ \tilde{x}_{21}\ \tilde{x}_{22}\ \tilde{x}_{31}\ \tilde{x}_{32}]$ of the transformed LP problem (9)-(10), the optimal solution to the original problem (7)-(8) can be retrieved via $[x_{11}\ x_{12}]' = \mathbf{Y}_1 \cdot [\tilde{x}_{11}\ \tilde{x}_{12}]'$, $[x_{21}\ x_{22}]' = \mathbf{Y}_2 \cdot [\tilde{x}_{21}\ \tilde{x}_{22}]'$, and $[x_{31}\ x_{32}]' = \mathbf{Y}_3 \cdot [\tilde{x}_{31}\ \tilde{x}_{32}]'$ by the three entities individually.

### B.2 Horizontally Partitioned LP Problem

For the LP problem (7)-(8), assuming $a_{11}$-$a_{14}$ and $b_{11}$-$b_{12}$ are owned by Entity I, $a_{21}$-$a_{24}$ and $b_{21}$-$b_{22}$ are owned by Entity II, while $a_{15}$-$a_{18}$, $a_{25}$-$a_{28}$, $a_{31}$-$a_{34}$, and $b_{31}$-$b_{32}$ are owned by Entity III. Each entity is unwilling to disclose the corresponding data. This is called a horizontally partitioned LP problem, as parameters in each row of constraints may be owned by different entities.

First, by enlarging the dimension of variables with additional slack variables $s_1$, $s_2$, and $s_3$, Entities I-III can equivalently transfer their corresponding inequality constraints into equalities. Thus, the original LP problem (7)-(8) is transformed into an equivalent LP problem (11)-(13). Note that coefficients of slack variables $r_1$, $r_2$, and $r_3$ in constraint (12) can be any random positive numbers to make (12) valid.

$$min\ [c_{11}\ c_{12}\ 0\ c_{21}\ c_{22}\ 0\ c_{31}\ c_{32}\ 0] \cdot$$
$$[x_{11}\ x_{12}\ s_1\ x_{21}\ x_{22}\ s_2\ x_{31}\ x_{32}\ s_3]' \quad (11)$$

$$\begin{bmatrix} a_{11} & a_{12} & 0 & 0 & 0 & 0 & 0 & 0 & 0 \\ a_{13} & a_{14} & r_1 & 0 & 0 & 0 & 0 & 0 & 0 \\ 0 & 0 & 0 & a_{21} & a_{22} & 0 & 0 & 0 & 0 \\ 0 & 0 & 0 & a_{23} & a_{24} & r_2 & 0 & 0 & 0 \\ a_{15} & a_{16} & 0 & a_{25} & a_{26} & 0 & a_{31} & a_{32} & 0 \\ a_{17} & a_{18} & 0 & a_{27} & a_{28} & 0 & a_{33} & a_{34} & r_3 \end{bmatrix} \cdot \begin{bmatrix} x_{11} \\ x_{12} \\ s_1 \\ x_{21} \\ x_{22} \\ s_2 \\ x_{31} \\ x_{32} \\ s_3 \end{bmatrix} \begin{matrix} = b_{11} \\ = b_{12} \\ = b_{21} \\ = b_{22} \\ = b_{31} \\ = b_{32} \end{matrix} \quad (12)$$

$$s_1, s_2, s_3 \ge 0 \quad (13)$$

Entities I-III can generate *2x2* random matrices $\mathbf{X}_1$, $\mathbf{X}_2$, and $\mathbf{X}_3$, respectively. Here "2" represents the number of constraints owned by each entity. That is, $\mathbf{X}_i = \begin{bmatrix} x_{i1} & x_{i2} \\ x_{i3} & x_{i4} \end{bmatrix}$ for $i = 1,2,3$. Thus, instead of disclosing $a_{11}$-$a_{14}$ and $b_{11}$-$b_{12}$, Entity I will provide $\mathbf{X}_1 \cdot \begin{bmatrix} b_{11} \\ b_{12} \end{bmatrix}$ and $\mathbf{X}_1 \cdot \begin{bmatrix} a_{11} & a_{12} & 0 & 0 & 0 & 0 & 0 & 0 \\ a_{13} & a_{14} & r_1 & 0 & 0 & 0 & 0 & 0 \end{bmatrix}$. Similarly, Entity II will provide $\mathbf{X}_2 \cdot \begin{bmatrix} b_{21} \\ b_{22} \end{bmatrix}$ and $\mathbf{X}_2 \cdot \begin{bmatrix} 0 & 0 & 0 & a_{21} & a_{22} & 0 & 0 & 0 \\ 0 & 0 & 0 & a_{23} & a_{24} & r_2 & 0 & 0 \end{bmatrix}$, instead of $a_{21}$-$a_{24}$ and $b_{21}$-$b_{22}$, while Entity III will provide $\mathbf{X}_3 \cdot \begin{bmatrix} b_{31} \\ b_{32} \end{bmatrix}$ and $\mathbf{X}_3 \cdot \begin{bmatrix} a_{15} & a_{16} & 0 & a_{25} & a_{26} & 0 & a_{31} & a_{32} & 0 \\ a_{17} & a_{18} & 0 & a_{27} & a_{28} & 0 & a_{33} & a_{34} & r_3 \end{bmatrix}$, instead of $a_{15}$-$a_{18}$, $a_{25}$-$a_{28}$, $a_{31}$-$a_{34}$, and $b_{31}$-$b_{32}$.

With the encrypted information, the original LP problem can be formulated as (14)-(16).

$$min\ [c_{11}\ c_{12}\ 0\ c_{21}\ c_{22}\ 0\ c_{31}\ c_{32}\ 0] \cdot$$
$$[x_{11}\ x_{12}\ s_1\ x_{21}\ x_{22}\ s_2\ x_{31}\ x_{32}\ s_3]' \quad (14)$$

$$\begin{bmatrix} \mathbf{X}_1 & 0 & 0 \\ 0 & \mathbf{X}_2 & 0 \\ 0 & 0 & \mathbf{X}_3 \end{bmatrix} \cdot \begin{bmatrix} a_{11} & a_{12} & 0 & 0 & 0 & 0 & 0 & 0 & 0 \\ a_{13} & a_{14} & r_1 & 0 & 0 & 0 & 0 & 0 & 0 \\ 0 & 0 & 0 & a_{21} & a_{22} & 0 & 0 & 0 & 0 \\ 0 & 0 & 0 & a_{23} & a_{24} & r_2 & 0 & 0 & 0 \\ a_{15} & a_{16} & 0 & a_{25} & a_{26} & 0 & a_{31} & a_{32} & 0 \\ a_{17} & a_{18} & 0 & a_{27} & a_{28} & 0 & a_{33} & a_{34} & r_3 \end{bmatrix} \cdot \begin{bmatrix} x_{11} \\ x_{12} \\ s_1 \\ x_{21} \\ x_{22} \\ s_2 \\ x_{31} \\ x_{32} \\ s_3 \end{bmatrix} =$$

$$\begin{bmatrix} \mathbf{X}_1 & 0 & 0 \\ 0 & \mathbf{X}_2 & 0 \\ 0 & 0 & \mathbf{X}_3 \end{bmatrix} \cdot [b_{11}\ b_{12}\ b_{21}\ b_{22}\ b_{31}\ b_{32}]' \quad (15)$$

$$s_1, s_2, s_3 \ge 0 \quad (16)$$

Comparing (7)-(8) with (14)-(16), it is clear that the extrema are equal, and the optimal solution to the original LP problem (7)-(8) is the same as that of transformed problem (14)-(16).

### C. Privacy Preserving ED Problem

The ED problem (1)-(6) can be written in the matrix form (17)-(19).

$$max\ [-\mathbf{c}_1\ \cdots\ -\mathbf{c}_G\ \mathbf{d}_1\ \cdots\ \mathbf{d}_D\ \mathbf{0}] \cdot$$
$$[\mathbf{P}_1\ \cdots\ \mathbf{P}_G\ \mathbf{D}_1\ \cdots\ \mathbf{D}_D\ \mathbf{\theta}]' \quad (17)$$

$$\begin{bmatrix} \mathbf{E}_1 & \cdots & \mathbf{0} & \mathbf{0} & \cdots & \mathbf{0} & \mathbf{0} \\ & \ddots & & & \ddots & & \vdots \\ \mathbf{0} & \cdots & \mathbf{E}_G & \mathbf{0} & \cdots & \mathbf{0} & \mathbf{0} \\ \mathbf{0} & \cdots & \mathbf{0} & \mathbf{F}_1 & \cdots & \mathbf{0} & \mathbf{0} \\ & \ddots & & & \ddots & & \vdots \\ \mathbf{0} & \cdots & \mathbf{0} & \mathbf{0} & \cdots & \mathbf{F}_D & \mathbf{0} \\ \mathbf{0} & \cdots & \mathbf{0} & \mathbf{0} & \cdots & \mathbf{0} & \mathbf{G} \cdot \mathbf{KL} \\ \mathbf{0} & \cdots & \mathbf{0} & \mathbf{0} & \cdots & \mathbf{0} & -\mathbf{G} \cdot \mathbf{KL} \end{bmatrix} \cdot \begin{bmatrix} \mathbf{P}_1 \\ \vdots \\ \mathbf{P}_G \\ \mathbf{D}_1 \\ \vdots \\ \mathbf{D}_D \\ \boldsymbol{\theta} \end{bmatrix} \leq \begin{bmatrix} \mathbf{M}_1 \\ \vdots \\ \mathbf{M}_G \\ \mathbf{N}_1 \\ \vdots \\ \mathbf{N}_D \\ \overline{\mathbf{PL}} \\ \overline{\mathbf{PL}} \end{bmatrix} \quad (18)$$

$$[\mathbf{KP}_1 \cdots \mathbf{KP}_G \; -\mathbf{KD}_1 \cdots -\mathbf{KD}_D \; -\mathbf{B}] \cdot \begin{bmatrix} \mathbf{P}_1 \\ \vdots \\ \mathbf{P}_G \\ \mathbf{D}_1 \\ \vdots \\ \mathbf{D}_D \\ \boldsymbol{\theta} \end{bmatrix} = \mathbf{0} \quad (19)$$

In the ED problem (17)-(19), each GENCO $i$ owns $\mathbf{c}_i$, $\mathbf{E}_i$, $\mathbf{M}_i$, and $\mathbf{KP}_i$ with decision variables $\mathbf{P}_i$, each LSE $j$ owns $\mathbf{d}_j$, $\mathbf{F}_j$, $\mathbf{N}_j$, and $\mathbf{KD}_j$ with decision variables $\mathbf{D}_j$, and the ISO owns $(\mathbf{G} \cdot \mathbf{KL})$, $\mathbf{B}$, and $\overline{\mathbf{PL}}$ with decision variables $\boldsymbol{\theta}$. Thus, the hybrid vertical and horizontal partition is applied to achieve the privacy-preserving ED model.

*1: Apply the Vertical Partition on the Original ED Problem (17)-(19).*

1.1 Each GENCO $i$ provides its encrypted matrices $(\mathbf{E}_i \cdot \mathbf{YG}_i)$, $(\mathbf{KP}_i \cdot \mathbf{YG}_i)$, and $(\mathbf{c}_i \cdot \mathbf{YG}_i)$, by generating its own $n_i \times n_i$ random positive matrix $\mathbf{YG}_i$. Each LSE $j$ provides its encrypted matrices $(\mathbf{F}_j \cdot \mathbf{YD}_j)$, $(\mathbf{KD}_j \cdot \mathbf{YD}_j)$, and $(\mathbf{d}_j \cdot \mathbf{YD}_j)$, by generating its own $n_j \times n_j$ random positive matrix $\mathbf{YD}_j$. The ISO provides its encrypted matrices $(\mathbf{G} \cdot \mathbf{KL} \cdot \mathbf{Y}_\theta)$ and $(\mathbf{B} \cdot \mathbf{Y}_\theta)$, by generating its own $(T \cdot (B-1)) \times (T \cdot (B-1))$ random positive matrix $\mathbf{Y}_\theta$.

1.2 With all information provided in 1.1, the transformed ED problem can be formulated as (20)-(22), which is a vertical partition privacy-preserving formulation of (17)-(19). That is, cost coefficients and left-hand-side coefficients of all constraints are encrypted via random numbers privately owned by individual entities.

$$\max \; [(-\mathbf{c}_1 \cdot \mathbf{YG}_1) \; \cdots \; (-\mathbf{c}_G \cdot \mathbf{YG}_G) \; (\mathbf{d}_1 \cdot \mathbf{YD}_1) \; \cdots \; (\mathbf{d}_D \cdot \mathbf{YD}_D) \; \mathbf{0}] \cdot [\widetilde{\mathbf{P}}_1 \cdots \widetilde{\mathbf{P}}_G \; \widetilde{\mathbf{D}}_1 \cdots \widetilde{\mathbf{D}}_D \; \widetilde{\boldsymbol{\theta}}]' \quad (20)$$

$$\begin{bmatrix} \mathbf{E}_1 \cdot \mathbf{YG}_1 & \cdots & \mathbf{0} & \mathbf{0} & \cdots & \mathbf{0} & \mathbf{0} \\ & \ddots & & & \ddots & & \vdots \\ \mathbf{0} & \cdots & \mathbf{E}_G \cdot \mathbf{YG}_G & \mathbf{0} & \cdots & \mathbf{0} & \mathbf{0} \\ \mathbf{0} & \cdots & \mathbf{0} & \mathbf{F}_1 \cdot \mathbf{YD}_1 & \cdots & \mathbf{0} & \mathbf{0} \\ & \ddots & & & \ddots & & \vdots \\ \mathbf{0} & \cdots & \mathbf{0} & \mathbf{0} & \cdots & \mathbf{F}_D \cdot \mathbf{YD}_D & \mathbf{0} \\ \mathbf{0} & \cdots & \mathbf{0} & \mathbf{0} & \cdots & \mathbf{0} & \mathbf{G} \cdot \mathbf{KL} \cdot \mathbf{Y}_\theta \\ \mathbf{0} & \cdots & \mathbf{0} & \mathbf{0} & \cdots & \mathbf{0} & -\mathbf{G} \cdot \mathbf{KL} \cdot \mathbf{Y}_\theta \end{bmatrix} \cdot \begin{bmatrix} \widetilde{\mathbf{P}}_1 \\ \vdots \\ \widetilde{\mathbf{P}}_G \\ \widetilde{\mathbf{D}}_1 \\ \vdots \\ \widetilde{\mathbf{D}}_D \\ \widetilde{\boldsymbol{\theta}} \end{bmatrix} \leq \begin{bmatrix} \mathbf{M}_1 \\ \vdots \\ \mathbf{M}_G \\ \mathbf{N}_1 \\ \vdots \\ \mathbf{N}_D \\ \overline{\mathbf{PL}} \\ \overline{\mathbf{PL}} \end{bmatrix} \quad (21)$$

$$[(\mathbf{KP}_1 \cdot \mathbf{YG}_1) \; \cdots \; (\mathbf{KP}_G \cdot \mathbf{YG}_G) \; -(\mathbf{KD}_1 \cdot \mathbf{YD}_1) \; \cdots \; -(\mathbf{KD}_D \cdot \mathbf{YD}_D) \; -\mathbf{B} \cdot \mathbf{Y}_\theta] \cdot [\widetilde{\mathbf{P}}_1 \cdots \widetilde{\mathbf{P}}_G \; \widetilde{\mathbf{D}}_1 \cdots \widetilde{\mathbf{D}}_D \; \widetilde{\boldsymbol{\theta}}]' = \mathbf{0} \quad (22)$$

*2: Apply the Horizontal Partition on the Transformed ED Problem (20)-(22).*

2.1 Each GENCO $i$ introduces $m_i$ slack variables $\widetilde{\mathbf{sG}}_i$ to convert its $m_i$ inequality constraints into equality, then provides its encrypted matrices $(\mathbf{XG}_i \cdot \mathbf{E}_i \cdot \mathbf{YG}_i)$, $(\mathbf{XG}_i \cdot \mathbf{RG}_i)$, and $(\mathbf{XG}_i \cdot \mathbf{M}_i)$, by generating one $m_i \times m_i$ random matrix $\mathbf{XG}_i$ and one $m_i \times m_i$ random positive diagonal matrix $\mathbf{RG}_i$. Each LSE $j$ introduces $m_j$ slack variables $\widetilde{\mathbf{sD}}_j$ to convert its $m_j$ inequality constraints into equality, then provides its encrypted matrices $(\mathbf{XD}_j \cdot \mathbf{F}_j \cdot \mathbf{YD}_j)$, $(\mathbf{XD}_j \cdot \mathbf{RD}_j)$, and $(\mathbf{XD}_j \cdot \mathbf{N}_j)$, by generating one $m_j \times m_j$ random matrix $\mathbf{XD}_j$ and one $m_j \times m_j$ random positive diagonal matrix $\mathbf{RD}_j$. The ISO introduces $(2T \cdot L)$ slack variables $\widetilde{\mathbf{sL}}_1$ and $\widetilde{\mathbf{sL}}_2$ to convert its $(2T \cdot L)$ inequality constraints into equality, and provides its encrypted matrices $(\mathbf{X}_{l1} \cdot \mathbf{G} \cdot \mathbf{KL} \cdot \mathbf{Y}_\theta)$, $(\mathbf{X}_{l2} \cdot \mathbf{G} \cdot \mathbf{KL} \cdot \mathbf{Y}_\theta)$, $(\mathbf{X}_{l1} \cdot \mathbf{R}_{l1})$, $(\mathbf{X}_{l2} \cdot \mathbf{R}_{l2})$, $(\mathbf{X}_{l1} \cdot \overline{\mathbf{PL}})$, and $(\mathbf{X}_{l2} \cdot \overline{\mathbf{PL}})$, as well as further encrypted matrices $(\mathbf{X}_{bl} \cdot \mathbf{KP}_i \cdot \mathbf{YG}_i)$, $(\mathbf{X}_{bl} \cdot \mathbf{KD}_j \cdot \mathbf{YD}_j)$, and $(\mathbf{X}_{bl} \cdot \mathbf{B} \cdot \mathbf{Y}_\theta)$, by generating two $(T \cdot L) \times (T \cdot L)$ random positive matrices $\mathbf{X}_{l1}$ and $\mathbf{X}_{l2}$, one $(T \cdot B) \times (T \cdot B)$ random positive matrix $\mathbf{X}_b$, and two $(T \cdot L) \times (T \cdot L)$ random positive diagonal matrices $\mathbf{R}_{l1}$ and $\mathbf{R}_{l2}$.

2.2 With all information provided in 2.1, the transformed ED problem can be formulated as (23)-(25), which is a horizontal partition privacy-preserving formulation of (20)-(22). That is, coefficients and right-hand-sides of all constraints are further encrypted via random numbers. Note that when converting an inequality constraint into equality, adding a non-negative slack variable or the product of a positive random number and the non-negative slack variable are mathematically equivalent. Thus, $\mathbf{RG}_i$, $\mathbf{RD}_j$, $\mathbf{R}_{l1}$, and $\mathbf{R}_{l2}$ are used to further encrypt each constraint set.

$$\max[(-\mathbf{c}_1 \cdot \mathbf{YG}_1) \cdots (-\mathbf{c}_G \cdot \mathbf{YG}_G) \; (\mathbf{d}_1 \cdot \mathbf{YD}_1) \cdots (\mathbf{d}_D \cdot \mathbf{YD}_D) \; \mathbf{0}] \cdot [\widetilde{\mathbf{P}}_1 \cdots \widetilde{\mathbf{P}}_G \; \widetilde{\mathbf{D}}_1 \cdots \widetilde{\mathbf{D}}_D \; \widetilde{\boldsymbol{\theta}} \; \widetilde{\mathbf{sG}}_1 \cdots \widetilde{\mathbf{sG}}_G \; \widetilde{\mathbf{sD}}_1 \cdots \widetilde{\mathbf{sD}}_D \; \widetilde{\mathbf{sL}}_1 \; \widetilde{\mathbf{sL}}_2]' \quad (23)$$

$$\widetilde{\mathbf{sG}}_1, \cdots, \widetilde{\mathbf{sG}}_G, \; \widetilde{\mathbf{sD}}_1, \cdots, \widetilde{\mathbf{sD}}_D, \; \widetilde{\mathbf{sL}}_1, \; \widetilde{\mathbf{sL}}_2 \geq 0 \quad (25)$$

*3: Solve the Transformed Problem (23)-(25) and Reconstruct the Optimal Solution to the Original ED Problem (17)-(19).*

3.1 Solve the transformed ED problem (23)-(25). Its optimal solution is $[\widetilde{\mathbf{P}}_1^* \cdots \widetilde{\mathbf{P}}_G^* \; \widetilde{\mathbf{D}}_1^* \cdots \widetilde{\mathbf{D}}_D^* \; \widetilde{\boldsymbol{\theta}}^* \; \widetilde{\mathbf{sG}}_1^* \cdots \widetilde{\mathbf{sG}}_G^* \; \widetilde{\mathbf{sD}}_1^* \cdots \widetilde{\mathbf{sD}}_D^* \; \widetilde{\mathbf{sL}}_1^* \; \widetilde{\mathbf{sL}}_2^*]$ and dual variable solution is $[\widetilde{\boldsymbol{\alpha}}_1^* \cdots \widetilde{\boldsymbol{\alpha}}_G^* \; \widetilde{\boldsymbol{\beta}}_1^* \cdots \widetilde{\boldsymbol{\beta}}_D^* \; \widetilde{\boldsymbol{\mu}}_1^* \; \widetilde{\boldsymbol{\mu}}_2^* \; \widetilde{\boldsymbol{\lambda}}^*]$, where $\widetilde{\boldsymbol{\alpha}}_i^*$ is an $m_i \times 1$

$$\begin{bmatrix} \mathbf{XG}_1 \cdot \mathbf{E}_1 \cdot \mathbf{YG}_1 & \cdots & \mathbf{0} & \mathbf{0} & \cdots & \mathbf{0} & \mathbf{0} & \mathbf{XG}_1 \cdot \mathbf{RG}_1 & \cdots & \mathbf{0} & \mathbf{0} & \cdots & \mathbf{0} & \mathbf{0} & \mathbf{0} \\ & \ddots & & & \ddots & & \vdots & & \ddots & & & \ddots & & \vdots & \vdots \\ \mathbf{0} & \cdots & \mathbf{XG}_G \cdot \mathbf{E}_G \cdot \mathbf{YG}_G & \mathbf{0} & \cdots & \mathbf{0} & \mathbf{0} & \mathbf{0} & \cdots & \mathbf{XG}_G \cdot \mathbf{RG}_G & \mathbf{0} & \cdots & \mathbf{0} & \mathbf{0} & \mathbf{0} \\ \mathbf{0} & \cdots & \mathbf{0} & \mathbf{XD}_1 \cdot \mathbf{F}_1 \cdot \mathbf{YD}_1 & \cdots & \mathbf{0} & \mathbf{0} & \mathbf{0} & \cdots & \mathbf{0} & \mathbf{XD}_1 \cdot \mathbf{RD}_1 & \cdots & \mathbf{0} & \mathbf{0} & \mathbf{0} \\ & \ddots & & & \ddots & & \vdots & & \ddots & & & \ddots & & \vdots & \vdots \\ \mathbf{0} & \cdots & \mathbf{0} & \mathbf{0} & \cdots & \mathbf{XD}_D \cdot \mathbf{F}_D \cdot \mathbf{YD}_D & \mathbf{0} & \mathbf{0} & \cdots & \mathbf{0} & \mathbf{0} & \cdots & \mathbf{XD}_D \cdot \mathbf{RD}_D & \mathbf{0} & \mathbf{0} \\ \mathbf{0} & \cdots & \mathbf{0} & \mathbf{0} & \cdots & \mathbf{0} & \mathbf{X}_{l1} \cdot \mathbf{G} \cdot \mathbf{KL} \cdot \mathbf{Y}_\theta & \mathbf{0} & \cdots & \mathbf{0} & \mathbf{0} & \cdots & \mathbf{0} & \mathbf{X}_{l1} \cdot \mathbf{R}_{l1} & \mathbf{0} \\ \mathbf{0} & \cdots & \mathbf{0} & \mathbf{0} & \cdots & \mathbf{0} & -\mathbf{X}_{l2} \cdot \mathbf{KL} \cdot \mathbf{G} \cdot \mathbf{Y}_\theta & \mathbf{0} & \cdots & \mathbf{0} & \mathbf{0} & \cdots & \mathbf{0} & \mathbf{0} & \mathbf{X}_{l2} \cdot \mathbf{R}_{l2} \\ \mathbf{X}_{bl} \cdot \mathbf{KP}_1 \cdot \mathbf{YG}_1 & \cdots & \mathbf{X}_{bl} \cdot \mathbf{KP}_G \cdot \mathbf{YG}_G & -\mathbf{X}_{bl} \cdot \mathbf{KD}_1 \cdot \mathbf{YD}_1 & \cdots & -\mathbf{X}_{bl} \cdot \mathbf{KD}_D \cdot \mathbf{YD}_D & -\mathbf{X}_{bl} \cdot \mathbf{B} \cdot \mathbf{Y}_\theta & \mathbf{0} & \cdots & \mathbf{0} & \mathbf{0} & \cdots & \mathbf{0} & \mathbf{0} & \mathbf{0} \end{bmatrix}$$

$$[\widetilde{\mathbf{P}}_1 \cdots \widetilde{\mathbf{P}}_G \; \widetilde{\mathbf{D}}_1 \cdots \widetilde{\mathbf{D}}_D \; \widetilde{\boldsymbol{\theta}} \; \widetilde{\mathbf{sG}}_1 \cdots \widetilde{\mathbf{sG}}_G \; \widetilde{\mathbf{sD}}_1 \cdots \widetilde{\mathbf{sD}}_D \; \widetilde{\mathbf{sL}}_1 \; \widetilde{\mathbf{sL}}_2]' = [(\mathbf{XG}_1 \cdot \mathbf{M}_1) \; \cdots \; (\mathbf{XG}_G \cdot \mathbf{M}_G) \; (\mathbf{XD}_1 \cdot \mathbf{N}_1) \; (\mathbf{XD}_D \cdot \mathbf{N}_D) \; (\mathbf{X}_{l1} \cdot \overline{\mathbf{PL}}) \; (\mathbf{X}_{l2} \cdot \overline{\mathbf{PL}}) \; \mathbf{0}]' \quad (24)$$

vector, $\widetilde{\boldsymbol{\beta}}_j^*$ is an $m_j \times 1$ vector, $\widetilde{\boldsymbol{\mu}}_1^*$ and $\widetilde{\boldsymbol{\mu}}_2^*$ are $(T \cdot L) \times 1$ vectors, and $\widetilde{\boldsymbol{\lambda}}^*$ is an $(T \cdot B) \times 1$ vector.

3.2 Each GENCO $i$ calculates its optimal solution via $\mathbf{P}_i^* = \mathbf{YG}_i \cdot \widetilde{\mathbf{P}}_i^*$. Each LSE $j$ calculates its optimal solution via $\mathbf{D}_j^* = \mathbf{YD}_j \cdot \widetilde{\mathbf{D}}_j^*$. The ISO can calculate the final bus angle solution via $\boldsymbol{\theta}^* = \mathbf{Y}_\theta \cdot \widetilde{\boldsymbol{\theta}}^*$. The final LMP can be calculated via $\mathbf{LMP}^* = -\mathbf{X}_b{}' \cdot \widetilde{\boldsymbol{\lambda}}^*$.

The detailed procedure of the proposed privacy-preserving ED approach is described as follows:

*Step 1*: Each GENCO $i$ provides its encrypted matrices $(\mathbf{c}_i \cdot \mathbf{YG}_i)$, $(\mathbf{XG}_i \cdot \mathbf{E}_i \cdot \mathbf{YG}_i)$, $(\mathbf{XG}_i \cdot \mathbf{RG}_i)$, and $(\mathbf{XG}_i \cdot \mathbf{M}_i)$. Each LSE $j$ provides its encrypted matrices $(\mathbf{d}_j \cdot \mathbf{YD}_j)$, $(\mathbf{XD}_j \cdot \mathbf{F}_j \cdot \mathbf{YD}_j)$, $(\mathbf{XD}_j \cdot \mathbf{RD}_j)$, and $(\mathbf{XD}_j \cdot \mathbf{N}_j)$.

*Step 2*: The ISO provides its encrypted matrices $(\mathbf{X}_{l1} \cdot \mathbf{G} \cdot \mathbf{KL} \cdot \mathbf{Y}_\theta)$, $(\mathbf{X}_{l1} \cdot \mathbf{R}_{l1})$, $(\mathbf{X}_{l2} \cdot \mathbf{R}_{l2})$, $(\mathbf{X}_{l1} \cdot \overline{\mathbf{PL}})$, and $(\mathbf{X}_{l2} \cdot \overline{\mathbf{PL}})$, as well as further encrypted matrices $(\mathbf{X}_{bl} \cdot \mathbf{KP}_i \cdot \mathbf{YG}_i)$, $(\mathbf{X}_{bl} \cdot \mathbf{KD}_j \cdot \mathbf{YD}_j)$, and $(\mathbf{X}_{bl} \cdot \mathbf{B} \cdot \mathbf{Y}_\theta)$.

*Step 3*: The ISO or the third-party solves the transformed privacy-preserving ED problem (23)-(25). Its optimal solution is $[\widetilde{\mathbf{P}}_1^* \cdots \widetilde{\mathbf{P}}_G^* \; \widetilde{\mathbf{D}}_1^* \cdots \widetilde{\mathbf{D}}_D^* \; \widetilde{\boldsymbol{\theta}}^* \; \widetilde{sG}_1^* \cdots \widetilde{sG}_G^* \; \widetilde{sD}_1^* \cdots \widetilde{sD}_D^* \; \widetilde{sL}_1^* \; \widetilde{sL}_2^*]$ and dual variable solution is $[\widetilde{\boldsymbol{\alpha}}_1^* \cdots \widetilde{\boldsymbol{\alpha}}_G^* \; \widetilde{\boldsymbol{\beta}}_1^* \cdots \widetilde{\boldsymbol{\beta}}_D^* \; \widetilde{\boldsymbol{\mu}}_1^* \; \widetilde{\boldsymbol{\mu}}_2^* \; \widetilde{\boldsymbol{\lambda}}^*]$.

*Step 4*: The corresponding optimal solution to the transformed privacy-preserving ED problem (23)-(25) is sent back to individual GENCOs and LSEs to reconstruct the optimal solution to the original ED problem. That is, each GENCO $i$ can obtain its optimal generation dispatch $\mathbf{P}_i^* = \mathbf{YG}_i \cdot \widetilde{\mathbf{P}}_i^*$; each LSE $j$ can calculate its optimal load dispatch $\mathbf{D}_j^* = \mathbf{YD}_j \cdot \widetilde{\mathbf{D}}_j^*$; and the ISO can derive the optimal bus angle $\boldsymbol{\theta}^* = \mathbf{Y}_\theta \cdot \widetilde{\boldsymbol{\theta}}^*$ and the final LMP $\mathbf{LMP}^* = -\mathbf{X}_b{}' \cdot \widetilde{\boldsymbol{\lambda}}^*$.

In sum, the proposed privacy-preserving ED approach allows individual GENCOs and LSEs to mask their actual bidding information and physical data by multiplying with random numbers before submitting to ISOs/RTOs. This could avoid potential information leakage on financial data of market participants. In addition, the ISO can encrypt physical information on the transmission network topology, which would avoid potential information leakage on critical energy infrastructure, especially when the optimization problem is deployed under cloud computing environments. In the next section, numerical case studies will validate the effectiveness of the proposed approach, along with detailed analyses on its computation and communication costs.

### III. COMPUTATIONAL EXPERIMENT

The proposed privacy-preserving ED approach is tested on a 3-bus system and the IEEE 118-bus system to evaluate its effectiveness on encrypting critical energy infrastructure information and financial data of market participants while being able to reconstruct the optimal solution to the original ED problem.

#### A. *3-Bus System*

A 3-bus system is studied for a single hour, which includes 2 generators, 1 load, and 3 branches. Generator and load data, including their three-segment bidding information, are given in Table I. Transmission line data are given in Table II. U1 belongs to GENCO1, GENCO2 owns U2, and LSE1 serves load L1.

TABLE I GENERATOR/LOAD INFORMATION

| Entity | Unit/Load | Bus | Bidding Price ($/MWh) | Min Capacity (MW) | Max Capacity (MW) |
|---|---|---|---|---|---|
| GENCO1 | U1 | 1 | 10 | 10 | 90 |
|  |  |  | 15 | 0 | 90 |
|  |  |  | 18 | 0 | 90 |
| GENCO2 | U2 | 2 | 12 | 10 | 80 |
|  |  |  | 18 | 0 | 80 |
|  |  |  | 20 | 0 | 80 |
| LSE1 | L1 | 3 | 19 | 100 | 150 |
|  |  |  | 16 | 0 | 50 |
|  |  |  | 14 | 0 | 50 |

TABLE II TRANSMISSION INFORMATION

| Line | From | To | x (p.u.) | Capacity (MW) |
|---|---|---|---|---|
| 1 | 1 | 2 | 0.1 | 30 |
| 2 | 2 | 3 | 0.1 | 150 |
| 3 | 1 | 3 | 0.1 | 100 |

The single hour ED problem is formulated as (26), with bus 1 as the reference bus. The optimal solution of (26) is 110MW, 80MW, and 190MW for the two generators and the load. Bus angles of buses 2 and 3 and -1p.u. and -10p.u. Power flows on the three lines are 10MW, 90MW, and 100 MW. LMPs are 15$/MWh, 15.5$/MWh, and 16$/MWh for the three buses.

$$max \; (-10P_{11} - 15P_{12} - 18P_{13} - 12P_{21} - 18P_{22} - 20P_{23}$$
$$+ 19D_{11} + 16D_{12} + 14D_{32})$$
$$10 \leq P_{11} \leq 90, \quad 0 \leq P_{12} \leq 90, \quad 0 \leq P_{13} \leq 90$$
$$10 \leq P_{21} \leq 80, \quad 0 \leq P_{22} \leq 80, \quad 0 \leq P_{23} \leq 80$$
$$100 \leq D_{11} \leq 150, \quad 0 \leq D_{12} \leq 50, \quad 0 \leq D_{13} \leq 50$$
$$-30 \leq (0 - \theta_2)/0.1 \leq 30$$
$$-150 \leq (\theta_2 - \theta_3)/0.1 \leq 150$$
$$-100 \leq (0 - \theta_3)/0.1 \leq 100$$
$$P_{11} + P_{12} + P_{13} + 10\theta_2 + 10\theta_3 = 0$$
$$P_{21} + P_{22} + P_{23} - 20\theta_2 + 10\theta_3 = 0$$
$$-D_{11} - D_{12} - D_{13} + 10\theta_2 - 20\theta_3 = 0 \tag{26}$$

The proposed privacy-preserving ED formulation corresponding to (23)-(25) is shown in (27)-(28). Six slack variables are introduced to each GENCO/LSE for converting inequalities constraints into equalities, and six slack variables are introduced to convert power flow inequality constraints into equalities. The optimal solutions to (27)-(28) include 74.65, -58.16, 69.40, 163.25, -98.31, 37.91, 40.08, 401.18, and -99.44 for the nine masked variables of the two generators and the load, 33.03 and -47.84 for the two masked variables of bus angles, as well as 0, 800, 500, 125, 101.12, 0, 0, 111.11, 888.89, 0, 160, 0, 0, 116.28, 15.15, 363.64, 68.49, 0, 74.07, 250, 0, 181.82, 300, and 204.08 for the 24 masked slack variables. The optimal solutions to dual variables of the 27 constraints in (27)-(28) are 3.46, -6.59, -1.84, -0.37, -1.29, 10.31, 141.47, -173.32, 57, 171.17, 78.35, -234.96, 4.49, -3.90, -0.47, -2.72, 5.76, -3.49, 0.48, 6.40, -6.15, 0, 0, 0, -13.13, -16.22, and -0.95. Sending corresponding solutions back to individual entities, the final dispatch and LMP solutions can be derived by individual entities via the internal calculation with $\mathbf{YG}_i$, $\mathbf{YD}_j$, $\mathbf{Y}_\theta$, and $\mathbf{X}_b^T$, respectively. The retrieved final solutions are exactly the same as those obtained by the original ED problem (26). For instance, GENCO1 can retrieve its dispatch solution of 110MW by multiplying its own random matrix $\mathbf{YG}_1$ and the optimal solutions of $\tilde{P}_{11}$, $\tilde{P}_{12}$, and $\tilde{P}_{13}$ from (27)-(28) as shown in (29). In addition, the final bus angles and LMPs can be retrieved via (30)-(31) by the ISO with its own random matrices $\mathbf{Y}_\theta$ and $\mathbf{X}_b$.

$$\begin{bmatrix}P_{11}\\P_{12}\\P_{13}\end{bmatrix} = \mathbf{YG}_1 \cdot \begin{bmatrix}\tilde{P}_{11}\\\tilde{P}_{12}\\\tilde{P}_{13}\end{bmatrix} = \begin{bmatrix}0.38 & 0.05 & 0.93\\0.56 & 0.53 & 0.13\\0.07 & 0.77 & 0.57\end{bmatrix} \cdot \begin{bmatrix}74.65\\-58.16\\69.40\end{bmatrix} = \begin{bmatrix}90\\20\\0\end{bmatrix} \quad (29)$$

$$\begin{bmatrix}\theta_2\\\theta_3\end{bmatrix} = \mathbf{Y}_\theta \cdot \begin{bmatrix}\tilde{\theta}_2\\\tilde{\theta}_3\end{bmatrix} = \begin{bmatrix}0.94 & 0.67\\0.32 & 0.43\end{bmatrix} \cdot \begin{bmatrix}33.03\\-47.84\end{bmatrix} = \begin{bmatrix}-1\\-10\end{bmatrix} \quad (30)$$

$$\begin{bmatrix}LMP_1\\LMP_2\\LMP_3\end{bmatrix} = -\mathbf{X}_b' \cdot \begin{bmatrix}\text{dual variables}\\\text{of the last three}\\\text{constgraints in}(25)\end{bmatrix}$$
$$= -\begin{bmatrix}0.58 & 0.06 & 0.92\\0.44 & 0.87 & 0.22\\0.26 & 0.63 & 0.37\end{bmatrix}' \cdot \begin{bmatrix}-13.13\\-16.22\\-0.95\end{bmatrix} = \begin{bmatrix}15\\15.5\\16\end{bmatrix} \quad (31)$$

In order to perform the original ED calculation (26), each market participant has to submit its actual financial bidding prices and physical information to the ISO. For instance, GENCO1 will submit its three-segment bidding prices 10$/MWh, 15$/MWh, and 18 $/MWh, together with the three-segment dispatch ranges [10, 90]MW, [0, 90]MW, and [0, 90]MW. In comparison, under the proposed privacy-preserving ED approach (27)-(28), GENCO1 will provide $\mathbf{c}_1 \cdot \mathbf{YG}_1$ (32), $\mathbf{XG}_1 \cdot \mathbf{E}_1 \cdot \mathbf{YG}_1$ (33), $\mathbf{Kp}_1 \cdot \mathbf{YG}_1$ (34), $\mathbf{XG}_1 \cdot \mathbf{RG}_1$ (35), and $\mathbf{XG}_1 \cdot \mathbf{M}_1$ (36). In addition, after solving (27)-(28), GENCO1 will use the final optimal solution to calculate and public its final optimal dispatch via (29). Random matrices $\mathbf{XG}_1$, $\mathbf{YG}_1$, and $\mathbf{RG}_1$ are only used by GENCO1 for internal calculations, but not shared with others.

We further investigate if inferences may occur in the proposed privacy-preserving ED approach. That is, with released information on $\mathbf{c}_1 \cdot \mathbf{YG}_1$, $\mathbf{XG}_1 \cdot \mathbf{E}_1 \cdot \mathbf{YG}_1$, $\mathbf{Kp}_1 \cdot \mathbf{YG}_1$, $\mathbf{XG}_1 \cdot \mathbf{RG}_1$, $\mathbf{XG}_1 \cdot \mathbf{M}_1$, and $\mathbf{YG}_1 \cdot [74.65 \quad -58.16 \quad 69.40]'$ shown in (32)-(37), whether one can reveal actual values of $\mathbf{c}_1$ and $\mathbf{M}_1$. (34) and (37) provide six linear equations with nine random variables of $\mathbf{YG}_1$. Thus, one cannot retrieve exact values of $\mathbf{YG}_1$ and in turn cannot derive actual values of $\mathbf{c}_1$ from (32). In addition, (33) and (35) provide 66 bi-linear equations with 51 random variables from $\mathbf{XG}_1$, $\mathbf{RG}_1$, and $\mathbf{YG}_1$. Thus, one cannot retrieve exact values of $\mathbf{XG}_1$ and in turn cannot derive actual values of $\mathbf{M}_1$ from (36). Furthermore, the first encrypted constraint for GENCO1 is shown in (38), which cannot not be used to reveal the actual constraint of $P_{11} \leq 90$. It is clear that by providing (32)-(37), GENCO1 can successfully

$$max \; [-13.46 \quad -22.31 \quad -21.51 \quad -17.16 \quad -25.84 \quad -18.44 \quad 17.09 \quad 14.54 \quad 30.45 \quad 0 \quad 0]$$
$$\cdot [\tilde{P}_{11} \; \tilde{P}_{12} \; \tilde{P}_{13} \; \tilde{P}_{21} \; \tilde{P}_{22} \; \tilde{P}_{23} \; \tilde{D}_{11} \; \tilde{D}_{12} \; \tilde{D}_{13} \; \tilde{\theta}_2 \; \tilde{\theta}_3]'$$

$$\begin{bmatrix}\mathbf{A}_1 & 0 & 0 & 0 & \mathbf{A}_5 & 0 & 0 & 0 & 0\\0 & \mathbf{A}_2 & 0 & 0 & 0 & \mathbf{A}_6 & 0 & 0 & 0\\0 & 0 & \mathbf{A}_3 & 0 & 0 & 0 & \mathbf{A}_7 & 0 & 0\\0 & 0 & 0 & \mathbf{A}_4 & 0 & 0 & 0 & \mathbf{A}_8 & 0\\0 & 0 & 0 & & 0 & 0 & 0 & 0 & \mathbf{A}_9\\ & \mathbf{A}_{10} & & & 0 & 0 & 0 & 0 & 0\end{bmatrix} \cdot [\tilde{P}_{11} \; \tilde{P}_{12} \; \tilde{P}_{13} \; \tilde{P}_{21} \; \tilde{P}_{22} \; \tilde{P}_{23} \; \tilde{D}_{11} \; \tilde{D}_{12} \; \tilde{D}_{32} \; \tilde{\theta}_2 \; \tilde{\theta}_3 \; \widetilde{sG}_{11} \; \widetilde{sG}_{12} \; \widetilde{sG}_{13} \; \widetilde{sG}_{14} \; \widetilde{sG}_{15} \; \widetilde{sG}_{16} \; \widetilde{sG}_{21}$$

$\widetilde{sG}_{22} \; \widetilde{sG}_{23} \; \widetilde{sG}_{24} \; \widetilde{sG}_{25} \; \widetilde{sG}_{26} \; \widetilde{sD}_{11} \; \widetilde{sD}_{12} \; \widetilde{sD}_{13} \; \widetilde{sD}_{14} \; \widetilde{sD}_{15} \; \widetilde{sD}_{16} \; \widetilde{sL}_{11} \; \widetilde{sL}_{12} \; \widetilde{sL}_{13} \; \widetilde{sL}_{21} \; \widetilde{sL}_{22} \; \widetilde{sL}_{23}]' = [233 \; 121.7 \; 65 \; 74.8 \; 98.6 \; 69.9 \; 151.4 \; 86.2 \; 140.1 \; 111.8 \; 166.7 \; 197.4 \; 85.5 \; 156.5 \; 51.5 \; -17.5 \; 177 \; 105 \; 98.9 \; 195.7 \; 187.1 \; 265.1 \; 205.5 \; 157.7 \; 0 \; 0 \; 0]'$

$\widetilde{sG}_{11}, \widetilde{sG}_{12}, \widetilde{sG}_{13}, \widetilde{sG}_{14}, \widetilde{sG}_{15}, \widetilde{sG}_{16}, \widetilde{sG}_{21}, \widetilde{sG}_{22}, \widetilde{sG}_{23}, \widetilde{sG}_{24}, \widetilde{sG}_{25}, \widetilde{sG}_{26}, \widetilde{sD}_{11}, \widetilde{sD}_{12}, \widetilde{sD}_{13}, \widetilde{sD}_{14}, \widetilde{sD}_{15}, \widetilde{sD}_{16}, \widetilde{sL}_{11}, \widetilde{sL}_{12}, \widetilde{sL}_{13}, \widetilde{sL}_{21}, \widetilde{sL}_{22}, \widetilde{sL}_{23} \geq 0$

(27)

$$\mathbf{A}_1 = \begin{bmatrix}0.77 & 0.43 & 0.82\\0.16 & -0.08 & -0.47\\-0.02 & -0.14 & -0.57\\-0.54 & 0.01 & -0.44\\-0.08 & -0.58 & -0.42\\-0.03 & -0.50 & -0.46\end{bmatrix}, \mathbf{A}_2 = \begin{bmatrix}0.28 & 0.34 & 0.02\\-0.23 & -0.39 & 0.04\\0.10 & 0.02 & -0.03\\-0.02 & 0.29 & 0.24\\0.41 & 0.27 & 0.79\\0.49 & 0.82 & 0.42\end{bmatrix}, \mathbf{A}_3 = \begin{bmatrix}-0.13 & 0.03 & -0.05\\0.15 & 0.22 & 0.18\\-0.59 & -0.06 & -0.53\\-0.59 & -0.35 & -0.76\\-0.01 & 0.32 & 0.16\\0.65 & 0.34 & 1.06\end{bmatrix}, \mathbf{A}_4 = \begin{bmatrix}-6.39 & -5.57\\-0.61 & -3.60\\-0.50 & -3.13\\5.27 & 7.61\\-0.81 & 2.41\\5.81 & 6.44\end{bmatrix}$$

$$\mathbf{A}_5 = \begin{bmatrix}0.07 & 0.00 & 0.14 & 0.02 & 0.70 & 0.81\\0.02 & 0.06 & 0.08 & 0.01 & 0.47 & 0.93\\0.01 & 0.05 & 0.08 & 0.04 & 0.10 & 0.44\\0.01 & 0.08 & 0.02 & 0.11 & 0.60 & 0.21\\0.04 & 0.05 & 0.07 & 0.08 & 0.16 & 0.85\\0.04 & 0.06 & 0.03 & 0.02 & 0.12 & 0.79\end{bmatrix}, \mathbf{A}_6 = \begin{bmatrix}0.47 & 0.29 & 0.07 & 0.01 & 0.29 & 0.19\\0.04 & 0.01 & 0.01 & 0.07 & 0.47 & 0.24\\0.55 & 0.35 & 0.03 & 0.02 & 0.37 & 0.25\\0.05 & 0.21 & 0.05 & 0.04 & 0.43 & 0.13\\0.71 & 0.03 & 0.02 & 0.07 & 0.49 & 0.05\\0.65 & 0.42 & 0.08 & 0.01 & 0.40 & 0.10\end{bmatrix}, \mathbf{A}_8 = \begin{bmatrix}0.25 & 0.10 & 0.04\\0.05 & 0.16 & 0.41\\0.07 & 0.17 & 0.32\end{bmatrix}, \mathbf{A}_9 = \begin{bmatrix}0.19 & 0.77 & 0.93\\0.08 & 0.75 & 0.53\\0.15 & 0.37 & 0.67\end{bmatrix}$$

$$\mathbf{A}_7 = \begin{bmatrix}0.59 & 0.34 & 0.42 & 0.04 & 0.21 & 0.20\\0.69 & 0.15 & 0.22 & 0.06 & 0.55 & 0.18\\0.40 & 0.26 & 0.43 & 0.03 & 0.01 & 0.33\\0.08 & 0.32 & 0.50 & 0.07 & 0.04 & 0.30\\0.63 & 0.04 & 0.38 & 0.05 & 0.49 & 0.27\\0.26 & 0.06 & 0.49 & 0.02 & 0.44 & 0.02\end{bmatrix}, \mathbf{A}_{10} = \begin{bmatrix}0.59 & 0.78 & 0.95 & 0.07 & 0.08 & 0.06 & -1.08 & -0.78 & -1.83 & 9.13 & 4.09\\0.44 & 0.59 & 0.72 & 0.96 & 1.23 & 0.90 & -0.26 & -0.19 & -0.44 & -7.37 & -3.50\\0.26 & 0.35 & 0.42 & 0.69 & 0.89 & 0.65 & -0.43 & -0.31 & -0.74 & -5.44 & -3.58\end{bmatrix}$$

(28)

encrypt its actual financial bidding prices and dispatch ranges from being released. Furthermore, it is impossible for other entities or third-parties to reveal actual values of $\mathbf{c}_1$ and $\mathbf{M}_1$ with public information on (32)-(37).

$$\mathbf{c}_1 \cdot \mathbf{YG}_1 = [c_{11}\ c_{12}\ c_{13}] \cdot \begin{bmatrix} YG_{11} & YG_{12} & YG_{13} \\ YG_{14} & YG_{15} & YG_{16} \\ YG_{17} & YG_{18} & YG_{19} \end{bmatrix} = \begin{bmatrix} 13.46 \\ 22.31 \\ 21.51 \end{bmatrix}' \quad (32)$$

$$\mathbf{XG}_1 \cdot \mathbf{E}_1 \cdot \mathbf{YG}_1 = \begin{bmatrix} XG_{11} & XG_{12} & XG_{13} & XG_{14} & XG_{15} & XG_{16} \\ XG_{17} & XG_{18} & XG_{19} & XG_{110} & XG_{111} & XG_{112} \\ XG_{113} & XG_{114} & XG_{115} & XG_{116} & XG_{117} & XG_{118} \\ XG_{119} & XG_{120} & XG_{121} & XG_{122} & XG_{123} & XG_{124} \\ XG_{125} & XG_{126} & XG_{127} & XG_{128} & XG_{129} & XG_{130} \\ XG_{131} & XG_{132} & XG_{133} & XG_{134} & XG_{135} & XG_{136} \end{bmatrix} \cdot$$
$$\begin{bmatrix} 1 & 0 & 0 \\ -1 & 0 & 0 \\ 0 & 1 & 0 \\ 0 & -1 & 0 \\ 0 & 0 & 1 \\ 0 & 0 & -1 \end{bmatrix} \cdot \begin{bmatrix} YG_{11} & YG_{12} & YG_{13} \\ YG_{14} & YG_{15} & YG_{16} \\ YG_{17} & YG_{18} & YG_{19} \end{bmatrix} = \begin{bmatrix} 0.77 & 0.43 & 0.82 \\ 0.16 & -0.08 & -0.47 \\ -0.02 & -0.14 & -0.57 \\ -0.54 & 0.01 & -0.44 \\ -0.08 & -0.58 & -0.42 \\ -0.03 & -0.50 & -0.46 \end{bmatrix} \quad (33)$$

$$\mathbf{Kp}_1 \cdot \mathbf{YG}_1 = \begin{bmatrix} 1 & 1 & 1 \\ 0 & 0 & 0 \\ 0 & 0 & 0 \end{bmatrix} \cdot \begin{bmatrix} YG_{11} & YG_{12} & YG_{13} \\ YG_{14} & YG_{15} & YG_{16} \\ YG_{17} & YG_{18} & YG_{19} \end{bmatrix} = \begin{bmatrix} 1.01 & 1.35 & 1.63 \\ 0 & 0 & 0 \\ 0 & 0 & 0 \end{bmatrix} \quad (34)$$

$$\mathbf{XG}_1 \cdot \mathbf{RG}_1 = \begin{bmatrix} XG_{11} & XG_{12} & XG_{13} & XG_{14} & XG_{15} & XG_{16} \\ XG_{17} & XG_{18} & XG_{19} & XG_{110} & XG_{111} & XG_{112} \\ XG_{113} & XG_{114} & XG_{115} & XG_{116} & XG_{117} & XG_{118} \\ XG_{119} & XG_{120} & XG_{121} & XG_{122} & XG_{123} & XG_{124} \\ XG_{125} & XG_{126} & XG_{127} & XG_{128} & XG_{129} & XG_{130} \\ XG_{131} & XG_{132} & XG_{133} & XG_{134} & XG_{135} & XG_{136} \end{bmatrix} \cdot$$
$$\begin{bmatrix} RG_{11} & & & & & \\ & RG_{12} & & & & \\ & & RG_{13} & & & \\ & & & RG_{14} & & \\ & & & & RG_{15} & \\ & & & & & RG_{16} \end{bmatrix} = \begin{bmatrix} 0.07 & 0 & 0.14 & 0.02 & 0.70 & 0.81 \\ 0.02 & 0.06 & 0.08 & 0.01 & 0.47 & 0.93 \\ 0.01 & 0.05 & 0.08 & 0.04 & 0.10 & 0.44 \\ 0.01 & 0.08 & 0.02 & 0.11 & 0.60 & 0.21 \\ 0.04 & 0.05 & 0.07 & 0.08 & 0.16 & 0.85 \\ 0.04 & 0.06 & 0.03 & 0.02 & 0.12 & 0.79 \end{bmatrix} \quad (35)$$

$$\mathbf{XG}_1 \cdot \mathbf{M}_1 = \begin{bmatrix} XG_{11} & XG_{12} & XG_{13} & XG_{14} & XG_{15} & XG_{16} \\ XG_{17} & XG_{18} & XG_{19} & XG_{110} & XG_{111} & XG_{112} \\ XG_{113} & XG_{114} & XG_{115} & XG_{116} & XG_{117} & XG_{118} \\ XG_{119} & XG_{120} & XG_{121} & XG_{122} & XG_{123} & XG_{124} \\ XG_{125} & XG_{126} & XG_{127} & XG_{128} & XG_{129} & XG_{130} \\ XG_{131} & XG_{132} & XG_{133} & XG_{134} & XG_{135} & XG_{136} \end{bmatrix} \cdot \begin{bmatrix} M_{11} \\ M_{12} \\ M_{13} \\ M_{14} \\ M_{15} \\ M_{16} \end{bmatrix} =$$
$$[233\ \ 121.7\ \ 65\ \ 74.8\ \ 98.6\ \ 69.9]' \quad (36)$$

$$\mathbf{YG}_1 \begin{bmatrix} \tilde{P}_{11} \\ \tilde{P}_{12} \\ \tilde{P}_{13} \end{bmatrix} = \begin{bmatrix} YG_{11} & YG_{12} & YG_{13} \\ YG_{14} & YG_{15} & YG_{16} \\ YG_{17} & YG_{18} & YG_{19} \end{bmatrix} \cdot \begin{bmatrix} 74.65 \\ -58.16 \\ 69.40 \end{bmatrix} = \begin{bmatrix} 90 \\ 20 \\ 0 \end{bmatrix} \quad (37)$$

$0.77\tilde{P}_{11} + 0.43\tilde{P}_{12} + 0.82\tilde{P}_{13} + 0.07\widetilde{sG}_{11} + 0.001\widetilde{sG}_{12} + 0.14\widetilde{sG}_{13} + 0.02\widetilde{sG}_{14} + 0.70\widetilde{sG}_{15} + 0.81\widetilde{sG}_{16} = 233 \quad (38)$

B. *The IEEE 118-Bus System*

The IEEE 118-bus system is further tested for 24 hours while considering all prevailing ED constraints discussed in Section II, for illustrating the effectiveness of the proposed privacy-preserving ED approach for larger systems. The system includes 54 generating units and 91 loads. The following 5 cases are studied:

*Case 0*: The original ED, in which all entities provide their actual information.

*Case 1*: The proposed privacy-preserving ED, in which each GENCO/LSE owns only one generator/load.

*Cases 2-4*: The proposed privacy-preserving ED, in which each GENCO/LSE owns 2, 5, and 10 generators/loads, respectively. These three cases will illustrate how different sizes of GENCOs/LSEs may impact computation and communication costs of the proposed privacy-preserving ED.

Table III shows problem scales and computational performances of all 5 cases. All cases derive the same optimal solution of $\$2.0221*10^6$. The proposed privacy-preserving ED models in Cases 1-4 do not increase the total number of constraints as compared to Case 0, although inequality constraints in Case 0 are converted into equality constraints with additional slack variables in Cases 1-4. In turn, numbers of variables in Cases 1-4 are much higher than that of Case 0. In addition, constraint coefficient matrices in Cases 1-4 are much denser than that of Case 0. Furthermore, the larger the number of generators/loads each entity owns, the higher the density of the constraint coefficient matrix is. The computational performance in the last column of Table III shows that computing times of Cases 1-4 are at the similar level, which is about 12 to 16 times higher than that of Case 0.

TABLE III COMPUTATIONAL PERFORMANCE OF THE 118-BUS SYSTEM

| Case | # of Variables | # of Inequalities | # of Equalities | # of Nonzeros | CPU Time (s) |
|---|---|---|---|---|---|
| 0 | 11472 | 24372 | 5016 | 75432 | 0.649 |
| 1 | 35844 | 0 | 29388 | 1178364 | 8.266 |
| 2 | 35844 | 0 | 29388 | 1270164 | 8.657 |
| 3 | 35844 | 0 | 29388 | 1538924 | 8.766 |
| 4 | 35844 | 0 | 29388 | 1964724 | 10.594 |

As the proposed privacy-preserving ED approach in Cases 1-4 evolves generating random matrices for encrypting financial and physical information, multiple instances of Case 4 are performed to illustrate how random matrices may impact the computational performance. Fig. 1 shows computational times of 100 instances of Case 4. The shortest and the longest computing times among the 100 instances are 6.375 seconds and 13.406 seconds, respectively. The mean computing time of the 100 instances is 8.639 seconds, with the standard deviation of 1.562 seconds. It shows that the computing time is impacted by random matrices, while the computational performance may range from 10 to 20 times higher than that of Case 0.

The proposed privacy-preserving ED approach, by submitting encrypted matrices instead of actual financial and physical information, may also increase the communication cost. Table IV shows communication requirements of the 5 cases. In comparison, information exchange between entities and the ISO in Cases 1-4 is much higher than that of Case 0, and is more significant with a larger number of generators/loads owned by each entity. For instance, in Case 0, the number of data sent from all entities to ISO (including bidding prices and physical limitations) is 24,216, and the number of data sent from ISO to all entities (including market clear quantities and LMPs) is 11,496. On the other hand, in Case 1, the number of data sent from all entities to ISO is increased to 232,536, which is about 8.6 times higher than that of Case 0. The last two columns of Table IV show the size of exchanged data in all 5 cases. Fortunately, with a typical 10 Mbps bandwidth communication infrastructure [17] between market participants and the ISO, all these communication tasks can be done within seconds.

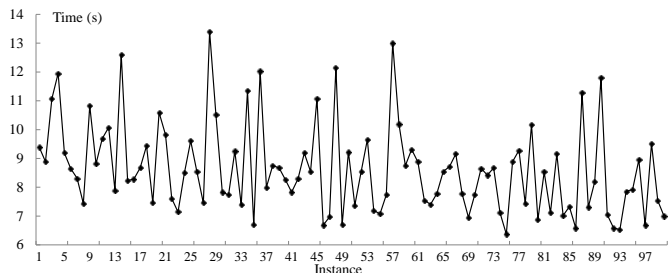

Fig. 1 Computational time of 100 instances for Case 4

TABLE IV COMMUNICATION COST

| | # of Exchanged Data | | Size of Exchanged Data (Mb) | |
|---|---|---|---|---|
| Case | From All Entities to ISO | From ISO to All Entities | From All Entities to ISO | From ISO to All Entities |
| 0 | 24,216 | 11,496 | 0.77 | 0.37 |
| 1 | 232,536 | 11,496 | 7.44 | 0.37 |
| 2 | 431,256 | 11,496 | 13.81 | 0.37 |
| 3 | 1,013,016 | 11,496 | 32.42 | 0.37 |
| 4 | 1,934,616 | 11,496 | 61.91 | 0.37 |

In sum, the following observations can be made from the above case studies:

1) The proposed privacy-preserving ED approach can retrieve the same optimal solution as the original ED approach, while avoiding disclosing actual bidding information and physical data.

2) For a same system, when a large number of generators/loads is owned by a single entity rather than multiple entities, the computational time does not increase noticeably while the communication cost would significantly increase.

3) Although the proposed privacy-preserving ED approach requires more computational time and a higher communication cost as compared to the original ED model, it is still tolerable under the current electricity market practice (i.e., less than 10 seconds for solving the privacy-preserving ED problem and a couple of seconds for the communication).

IV. CONCLUSION

This paper privacy-preserving ED approach in competitive electricity market, in which individual GENCOs and LSEs can mask their actual bidding information and physical data by multiplying with random numbers before submitting to ISOs/RTOs. This would avoid potential information leakage of potential financial data of market participants. In addition, the ISO can encrypt physical information on the transmission network topology, which would avoid potential information leakage on critical energy infrastructure, especially when the optimization problem is deployed under cloud computing environments. The transformed privacy-preserving ED model is still an LP problem, with the same number of total constraints and a larger number of variables as compared to the original ED problem. The optimal solution to the original ED problem can be retrieved from the optimal solution of the proposed privacy-preserving ED approach. Numerical case studies show that although computation and communication costs of the proposed privacy-preserving ED approach are higher than the original ED, it is still tolerable under the current electricity market practice. The future work will investigate the application of the proposed approach on other power systems optimization problems such as unit commitment and long-term investment planning.